\documentclass[12pt]{article}
\usepackage[utf8]{inputenc}
\usepackage[T1]{fontenc}
\usepackage{amssymb}
\usepackage{amsmath}
\usepackage{xcolor}
\usepackage{graphicx}
\usepackage{float}
\usepackage{amsthm}
\usepackage{indentfirst}
\usepackage{fancybox}
\usepackage[top=2cm,bottom=2cm,left=2cm,right=2cm]{geometry}
\usepackage{hyperref}
\hypersetup{
    colorlinks=true,
    linkcolor=blue,
    filecolor=magenta,      
    urlcolor=cyan,
    pdftitle={Sharelatex Example},
    bookmarks=true,
    pdfpagemode=FullScreen,
    }

\usepackage{tikz}

\usetikzlibrary{arrows,snakes,automata,backgrounds,calc, positioning,quotes}
\usetikzlibrary{decorations.pathmorphing}
\usepackage{pst-eucl,pstricks-add,pst-fun,pst-tree}

{\theoremstyle{definition}
\newtheorem{theo}{Theorem}[section]}

{\theoremstyle{definition}
}

{\theoremstyle{definition}
}

{\theoremstyle{definition}
\newtheorem*{dem}{Proof}}

\usepackage{comment}

\title{Analysis of Least square estimator for simple Linear Regression with a uniform distribution error}

\author{M.R. Jlibene, \and S. Taoufik, \and S. Benjelloun
\\
{\small MSDA, Université Mohammed VI Polytechnique, Hay Moulay Rachid Benguerir, Morocco. } }
\date{December 2020}

\usepackage{natbib}
\usepackage{graphicx}

\begin{document}

\maketitle

\begin{abstract}
    We study the least square estimator, in the framework of simple linear regression, when the deviance term $\epsilon$ with respect to the linear model is modeled by a uniform distribution. In particular, we give the law of this estimator, and prove some convergence properties.  
\end{abstract}
\section{Introduction}

we consider the standard linear regression problem $Y = \beta_0 + \beta_1 X + \epsilon$ where $X$ is a  random or deterministic  variable in $\mathbb{R}$, $\beta \in \mathbb{R}$ the vector of model coefficients, and $\epsilon$ a scalar centered random vector, $\mathbb{E}(\epsilon) = 0$, modeling a random perturbation term.

Given a $n$-sample of data points we adopt the following matrix notation :
$$Y =  [1, x] \beta + \epsilon $$ 
\noindent where :
\begin{itemize}
\item $x$ is the data vector in $\mathbb{R}^n$. Furthermore we suppose that the vector $x$ is not co-linear with $[1, 1...,1]^T \in \mathbb{R}^n$. We note $X$ the matrix $[1, x]$.

\item $\beta = [\beta_0, \beta_1]^T \in \mathbb{R}^2$ are the coefficients vector of the regression model.

\item $\epsilon = (\epsilon_i)_i \in \mathbb{R}^n$ is the vector for the disturbance samples. We suppose that the $n$ identically distributed realisations are non-correlated, or eventually independent.
\end{itemize}

The Gauss-Markov theorem apply in our setting and the least square estimator (LSE) $\hat{\beta}_n = (X^T X)^{-1} X^T Y = \beta + (X^T X)^{-1} X^T \epsilon = M \epsilon $ is the Best Linear Unbiased Estimator (BLUE) for $\beta$. The least square estimator $\hat{\beta}_n$ can be also written as : 
\begin{center}
    $\hat{\beta}_n = \hat{X}PY = \hat{X}P(X\beta + \epsilon)$ \\
    $\hat{\beta}_n = \beta + \hat{X}P\epsilon$
\end{center}
\noindent with 
\begin{itemize}
    \item $P = X (X^T X)^{-1}X^T $ is the projection  in  $\mathbb{R}^n$ on $\text{Im}(X)$.
    \item $\hat{X}: \text{Im}(X) \subset \mathbb{R}^n \rightarrow \mathbb{R}^d$ is the inverse of $X$ restricted to $\text{Im}(X)$.
\end{itemize}

Hence, computing the distribution law of $\hat{\beta}_n$, boils down to the computation of the distribution law of $P\epsilon$. When the disturbance term $\epsilon$ has a centred Gaussian distribution $\mathcal{N}(0, \sigma^2)$, of known variance $\sigma^2$, with independent realisations making $(\epsilon_i)_i$ a Gaussian vector, the computation of this law is straightforward once the invert $(X^T X)^{-1}$ is computed : $\hat{\beta}_n \sim \mathcal{N}(\beta, (X^T X)^{-1} \sigma^2)$.

We have :
$$
X = \left[
\begin{array}{cc}
1 & x_1 \\
\vdots  & \vdots\\
1 & x_n 
\end{array}
\right]
\, , \,  X^T X = \left[
\begin{array}{cc}
n &  \sum x_i\\
\sum x_i &  \sum x_i^2
\end{array}
\right] 
$$
$$
\text{det}(X^T X) = n (\sum x_i^2) - (\sum x_i)^2 = \sum_{i \neq j} (x_j - x_i)^2
$$ 
$$
(X^T X)^{-1} = \frac{1}{\text{det}(X^T X)} \left[
\begin{array}{cc}
\sum x_i^2 &  -\sum x_i\\
-\sum x_i & n 
\end{array}
\right].
$$
$$
P_{i,j} = \frac{1}{\text{det}(X^T X)}\left[\sum x_k^2+ n \, x_i \, x_j - (x_j + x_i)\sum x_k\right]
$$
$$
\hat{X} =\frac{1}{det(X^TX)}\left[\begin{matrix}   \sum_{k\neq 1} x_k^2  - x_1 \sum_{k\neq 1} x_k & ... & \sum_{k\neq j} x_k^2  - x_j \sum_{k\neq j} x_k & ... &   \sum_{k\neq n} x_k^2  - x_n \sum_{k\neq n}  x_k \\
 n \, x_1   -  \sum  x_k & ... & n \, x_j - \sum  x_k & ...  & n \, x_n  -  \sum x_k 
\end{matrix}
\right]
$$

$$
\hat{\beta} = \beta + M \, \epsilon
$$
with 
{\small
$$
M = (X^T X)^{-1} X^T = \frac{1}{det(X^T X)} 
\left[\begin{matrix}   \sum_{k\neq 1} x_k^2  - x_1 \sum_{k\neq 1} x_k & ... & \sum_{k\neq j} x_k^2  - x_j \sum_{k\neq j} x_k & ... &   \sum_{k\neq n} x_k^2  - x_n \sum_{k\neq n}  x_k \\
 n \, x_1   -  \sum  x_k & ... & n \, x_j - \sum  x_k & ...  & n \, x_n  -  \sum x_k 
\end{matrix}
\right]
$$
}
$$
 = \frac{1}{det(X^T X)} 
\left[\begin{matrix}   \sum_{k\neq 1} x_k^2  - x_1 \sum_{k\neq 1} x_k & ... & \sum_{k\neq j} x_k^2  - x_j \sum_{k\neq j} x_k & ... &   \sum_{k\neq n} x_k^2  - x_n \sum_{k\neq n}  x_k \\
 n \, x_1   -  \sum  x_k & ... & n \, x_j - \sum  x_k & ...  & n \, x_n  -  \sum x_k 
\end{matrix}
\right]
$$

We note that 

$$
\sum_j M_{1,j}  = \frac{1}{det(X^T X)} \left[(n-1) \sum_i x_i ^2 -  \sum_{i\neq j} x_i x_j\right]
=\frac{1}{det(X^T X)} \sum_{i\neq j} (x_i-x_j)^2
= 1
$$

$$
\sum_j M_{2,j} = 0
$$


When $\sigma$ is unknown, from the maximum likelihood estimator one can construct the unbiased  estimator for $\sigma^2$ is given by $$\hat{\sigma}^2 =  \frac{1}{n-2} ||Y - \hat{Y} ||^2 $$
\noindent with $\hat{Y} = X\hat{\beta}$ .

We summarize here in the results for the Gaussian case :
\begin{itemize}
    \item Law for $\hat{\beta}$ : $\mathcal{N}\left( \beta, (X^T X)^{-1} \sigma^2 \right)$. In our case (simple linear regression, $d=1.5$) we have:
    \begin{itemize}
    \item $\hat{\beta}_{0}\sim \mathcal{N}\left(\beta_{0}, \frac{S_2}{nS_2-S_1^{2}}\sigma^2\right)$ where $S_1=\sum x_{i} $ and $S_2=\sum x_{i}^{2}$
    \item $\hat{\beta}_{1}\sim \mathcal{N}\left(\beta_{1}, \frac{n}{nS_2-S_1^{2}}\sigma^2\right)$
    
    \item $\operatorname{cov}(\hat{\beta}_{0}, \hat{\beta}_{1})=\frac{-S_1}{nS_2-S_1^{2}}\sigma^2$
    
    \item Hence, if we suppose $\sigma^2$ to be known, confidence intervals can be built using the normal laws:
$$\frac{\widehat{(\beta}_{0}-\beta_{0})\sqrt{nS_2-S_1^{2}}}{\sigma \sqrt{S_2}} \sim \mathcal{N}(0,1)$$

$$\frac{\widehat{(\beta}_{1}-\beta_{1})\sqrt{nS_2-S_1^{2}}}{\sigma \sqrt{n}} \sim \mathcal{N}(0,1)$$ 
    \item The (symetric) confidence intervals at level $1-\alpha$ for $\beta_0$ and $\beta_1$ are respectively :
  
    $$
    \left[\hat{\beta}_0 - q_{1- \frac{\alpha}{2}} \frac{\sigma \sqrt{S_2}}{\sqrt{nS_2-S_1^{2}}}, \hat{\beta}_0 + q_{1 - \frac{\alpha}{2}} \frac{\sigma \sqrt{S_2}}{\sqrt{nS_2-S_1^{2}}} \right]
    $$
    $$
    \left[ \hat{\beta}_1 - q_{ 1-  \frac{\alpha}{2}} \frac{\sigma \sqrt{n}}{\sqrt{nS_2-S_1^{2}}} , \hat{\beta}_1 + q_{1 - \frac{\alpha}{2}} \frac{\sigma \sqrt{n}}{\sqrt{nS_2-S_1^{2}}}  \right]
    $$
    where $q_{1 - \frac{\alpha}{2}}$ is the quantile at level $1 - \frac{\alpha}{2}$ for the normal law $\mathcal{N}(0,1)$.
    \item The statistical tests at significance level $\alpha$ for the hypothesis $H_0 : \beta_0=0$ and $H_0 : \beta_1=0$ have respectively the critical regions defined as :  
    $$
    | \hat{\beta}_1 |  > q_{1- \frac{\alpha}{2}} \frac{\sigma \sqrt{n}}{\sqrt{nS_2-S_1^{2}}}
    $$
    $$
    | \hat{\beta}_0 |  > q_{1 - \frac{\alpha}{2}}\frac{\sigma \sqrt{S_2}}{\sqrt{nS_2-S_1^{2}}}
    $$
    \end{itemize}
\end{itemize}

When the perturbation $\epsilon$ is not Gaussian, the above properties are still verified asymptotically (large sample size $n$).

we are interested in the case where the error is uniform, $\epsilon \sim U(-\theta, \theta)$, and when the sample is small. This may arise in many cases related to physical and experimental measures, where the sensibility of measurement equipment may be better modeled by random perturbations with a uniform law.

In the case where $\epsilon$ is uniform in $[-\theta,\theta]$, with independent realisations, the Gauss Markov theorem still apply in this case and the $LSE$ is the BLUE. We are interested by studying other properties of the LSE estimator in this case (law, confidence intervals, regression tests). In  the case where the error law parameters $\theta$ is known, a closed formula for the estimator distribution law can be derived by slightly generalising the Irwin-Hall distribution such in \hyperlink{sadoo}{Sadooghi, Nematollahi and Habibi (2007)}, or \hyperlink{bradl}{Bradley and Gupta (2002)}.

In section \ref{sec:LSE} we give this closed formula for the estimator law, in the case where the uniform law parameter are known $\epsilon \sim \mathcal{U}(-\theta, \theta)$ with known $\theta$. In section \ref{sec:example} we show a numerical test for illustration.


\section{Least Square estimator}\label{sec:LSE}

We start by introducing a generalisation of  the Erwin-Hall density function as follows : 

\begin{theo} Assume that $U_1$, $U_2$,...,$U_n$ are independents variables wih $U_k$ uniformly distributed on $[0, p_k]$ where $p_k> 0$, and let $U = \sum_{k=1}^n U_k$.
The probability density function for $U$ is given by
$$f_{U}(x;n)=\frac{1}{P_{n}(n-1) !} \sum_{k=0}^{n}(-1)^{k} \sum_{l=1}^{n_k}\left(x-S_{n k l}\right)_{+}^{n-1}$$
for $x \in [0,\sum_{k=1}^n p_k] $.\\
where $n_k=\left(\begin{array}{c}n \\k\end{array}\right)$ is the number of k-combinations, $l \in [0,n_k] $ is the index over the $n_k$ combinations, $(S_{n k l})_l$ are the sums over each $k-$combination
of  $\left(p_{1}, \ldots, p_{n}\right)$, and $P_{n}=\prod_{k=1}^n p_k$ and $x_+=max(0,x)$.
\end{theo}
\begin{dem}
See \hyperlink{sadoo}{Sadooghi, Nematollahi and Habibi (2007)}.$\square$
\end{dem}

The Erwin-Hall distribution corresponds to the special case $p_k =1$, so $S_{n k l}=k$ and 

$$
f(x;n)=\frac{1}{(n-1) !} \sum_{k=0}^{n}(-1)^{k} n_k \left( x-k \right)_{+}^{n-1}
$$

\textbf{Generalisation:}  let $V = \sum_{k=1}^n V_k$, where $V_1$, $V_2$...,$V_n$ are independents variables and $V_k\sim \mathrm{U}(-p_k,p_k)$ and independents. \\Let $U_k=\frac{V_k+p_k}{2}$ and $U=\sum_{k=1}^n U_k$. Clearly $U_k$ satisfies the Theorem 2.1 conditions. We have $V = \sum_{k=1}^n (2U_k - p_k)=2U -\sum_{k=1}^n p_k$  so the probability density function for $V$ is given by:
$$
f_V(x;n)=\frac{1}{2}f_U(\frac{x+\sum_{k=1}^n p_k}{2};n)=\frac{1}{2P_{n}(n-1) !} \sum_{k=0}^{n}(-1)^{k} \sum_{l=1}^{n_k}\left(\frac{x}{2}+\frac{(\sum_{i=1}^n p_i)}{2}-S_{n k l}\right)_{+}^{n-1}
$$

\textbf{Generalisation:} let  $W = P \epsilon = \sum_{k=1}^n p_k \epsilon_k$, $\epsilon_k\sim \mathrm{U}(-\theta,\theta).$

In this case $W = \sum_{k=1}^n W_k$, $W_k = p_k \epsilon_k \sim \mathrm{U}(-p_k \theta,p_k\theta)$, and

$$
f_W(x;n)=\frac{1}{2P_{n}\theta^{n}(n-1) !} \sum_{k=0}^{n}(-1)^{k} \sum_{l=1}^{n_k}\left(\frac{x}{2}+\theta.\frac{(\sum_{i=1}^n p_i)}{2}-\theta.S_{n k l}\right)_{+}^{n-1}
$$

In our case, as $\theta$ is known, we know that:
$\hat{\beta} = \beta + (X^T X)^{-1} X^T \epsilon.$

Hence, with : 
$$d=\text{det}(X^T X) = nS_2-S_1^{2},$$ 
$$p_{k}=S_2-x_{k}S_1,$$
$$p^{'}_{k}=x_{k}n-S_1,$$

\noindent we have

$$\hat{\beta}_0= \beta_0 + \frac{1}{d}\sum_{k=1}^n p_k \epsilon_k= \beta_0 + \frac{1}{d}W$$
$$\hat{\beta}_1= \beta_1 + \frac{1}{d}\sum_{k=1}^n p^{'}_k \epsilon_k= \beta_1 + \frac{1}{d}W'$$

$$f_{\hat{\beta}_{0}}(x;n)=d.f_W(d(x-{\beta}_{0});n)$$
$$f_{\hat{\beta}_{1}}(x;n)=d.f_{W'}(d(x-{\beta}_{1});n)$$
$$f_{\hat{\beta}_{0}}(x;n)=d.\frac{1}{P_{n}\theta^{n}(n-1) !} \sum_{k=0}^{n}(-1)^{k}\sum_{l=1}^{n_k}\left(d(x-{\beta}_{0})+\theta.d-\theta.S_{n k l}\right)_{+}^{n-1}$$

$$f_{\hat{\beta}_{1}}(x;n)=d.\frac{1}{P^{'}_{n}\theta^{n}(n-1) !} \sum_{k=0}^{n}(-1)^{k}\sum_{l=1}^{n_k}\left(d(x-{\beta}_{1})+-\theta.S^{'}_{n k l}\right)_{+}^{n-1}$$

\begin{theo}
If $\frac{\max _{1 \leq i \leq n}\left(|p_i|\right)}{\sqrt {d\,S_2}} \rightarrow 0$ when $n \rightarrow +\infty$, we have : 
$\sqrt{\frac{d}{S_2}}(\hat{\beta}_0-\beta_0)$ converges in law to the normal distribution $N(0, \frac{\theta^2}{3} )$. If $\frac{\max _{1 \leq i \leq n}\left(|p_i^{'}|\right)}{\sqrt {dn}} \rightarrow 0$ when $n \rightarrow +\infty$, we have $\sqrt{\frac{d}{n}}(\hat{\beta}_1-\beta_1)$  that converges in law to the normal distribution $N(0, \frac{\theta^2}{3} )$

\end{theo}

\begin{dem}     
We have $$\hat{\beta}_0-\beta_0=\frac{1}{d}\sum_{k=1}^n p_k \epsilon_k.$$

Hence
$\sqrt{\frac{d}{S_2}}(\hat{\beta}_0-\beta_0)=\frac{1}{\sqrt {dS_2}}\sum_{k=1}^n p_k \epsilon_k$.

The characteristic function for $\sqrt{\frac{d}{S_2}}(\hat{\beta}_0-\beta_0)$ is :
$$\phi_n(t)=\prod \frac{\sin( \frac{p_i\theta t}{\sqrt {dS_2}})}{\frac{p_i\theta t}{\sqrt {dS_2}}}$$

Let $t$ be a real number, as $n \rightarrow +\infty$ and given that $\frac{p_i\theta t}{\sqrt {dS_2}}\rightarrow 0$ we have 
$$\phi_n(t)=\prod (1- \frac{1}{6}(\frac{p_i\theta t}{\sqrt {dS_2}})^2+ o((\frac{p_i\theta t}{\sqrt {dS_2}})^2))$$
Hence
$$\ln(\phi_n(t))=\sum \ln(1- \frac{1}{6}(\frac{p_i\theta t}{\sqrt {dS_2}})^2+ o((\frac{p_i\theta t}{\sqrt {dS_2}})^2))$$

$$\ln(\phi_n(t))=\sum (- \frac{1}{6}(\frac{p_i\theta t}{\sqrt {dS_2}})^2+ o((\frac{p_i\theta t}{\sqrt {dS_2}})^2))$$

As $\sum p_i^2=dS_2$ we have $\sum o((\frac{p_i\theta t}{\sqrt {dS_2}})^2) \rightarrow 0$ and hence
$$\ln(\phi_n(t))=- \frac{1}{6}(\theta t)^2+ o(1)$$

We proved then that $\phi_n(t) \rightarrow \exp(- \frac{1}{6}(\theta t)^2)$ when $n \rightarrow +\infty$ which gives the result for $\beta_0$.

The same can be done for $\beta_1$ using $\sum p_i'^2=d\,n$. $\square$
\end{dem}

A similar result can be given for the limit law for the couple $\hat{\beta}$.

\begin{theo}
If $\frac{\max _{1 \leq i \leq n}\left(|p_i|\right)}{\sqrt {dS_2}} \rightarrow 0$ , $\frac{\max _{1 \leq i \leq n}\left(|p_i^{'}|\right)}{\sqrt {dn}} \rightarrow 0$  and $\frac{S_1}{nS_2} \rightarrow 0$ when $n \rightarrow +\infty$, then 
($\sqrt{\frac{d}{S_2}}(\hat{\beta}_0-\beta_0)$,$\sqrt{\frac{d}{n}}(\hat{\beta}_1-\beta_1)$) converges in law to $N_{2}\left(0, \frac{\theta^{2}}{3} \mathrm{I_2}\right)$ 
\end{theo}

\begin{dem}
The demonstration is similar to the previous theorem. It suffices to justify the convergence to $0$ for  $\sum (\frac{p_i \, p_i' }{d\sqrt {n \,S_2}})=-\frac{S_1}{n \, S_2}$ 
which is satisfied under the two conditions of the theorem. $\square$
\end{dem}

Let us now consider the question whether the conditions for the last theorems hold in practical situations.

If we consider the example where the $x_k$ are uniform over $[a,b]$ :
$$x_k=(b-a)\frac{k-1}{n-1}+a.$$ 
We can see that the two conditions of the previous theorems are indeed satisfied.

If the $x_k$ are iid random variables, uniformely distributed over an intervalle $[a,b]$, then using the strong law of large numbers we demontrate that the conditions for the two theorem are satisfied almost surely 
$$\frac{p_i}{\sqrt {dS_2}}=\frac{1}{\sqrt{n}} \cdot \frac{\frac{S_{2}}{n}-x_{k} \frac{S_{1}}{n}}{\sqrt{\left(\frac{S_{2}}{n}\right)^{2}-\left(\frac{S_{1}}{n}\right)^{2} \cdot \frac{S_{2}}{n}}},$$ 
and 
$$\frac{S_1}{nS_2}=\frac{1}{n}.\frac{S_1}{n}.\frac{1}{\frac{S_2}{n}}.$$

\subsection{Estimator for $\theta$}

We give an unbiased estimator for $\theta$.

\begin{theo} \label{th-theta} When $\theta$ is unknown, an unbiased estimator for $\theta^2$ is given by $$\hat{\theta}^2 =  \frac{3}{n-2} ||Y - \hat{Y} ||^2 $$
\noindent with $\hat{Y} = X\hat{\beta}=P Y$ 
\end{theo}

\begin{dem}
 we have $$||Y - \hat{Y} ||^2=(Y-\widehat{Y})^{T}(Y-\widehat{Y})=(Y-P Y)^{T}(Y-P Y)=Y^{T}\left(\mathbf{I}_{n}-P\right) Y$$
 
 $$
\begin{aligned}
\mathbb{E}(||Y - \hat{Y} ||^2) 
&=\mathbb{E}\left[\operatorname{trace}\left(Y^{T}\left(\mathbf{I}_{n}-P\right) Y\right)\right] \\
&=\mathbb{E}\left[\operatorname{trace}\left(\left(\mathbf{I}_{n}-P\right) Y Y^{T}\right)\right] \\
&=\operatorname{trace}\left[\left(\mathbf{I}_{n}-P\right) \mathbb{E}\left(Y Y^{T}\right)\right] \\
&=\operatorname{trace}\left[\left(\mathbf{I}_{n}-P\right)\left(\frac{\theta^2}{3} \mathbf{I}_{n}+X \beta \beta^{T} X^{T}\right)\right]\\
&=\operatorname{trace}\left[\frac{\theta^2}{3}\left(\mathbf{I}_{n}-P\right)\right] \\
&=\frac{\theta^2}{3}(n-2). \square
\end{aligned} 
$$

\end{dem}

\section{Numerical simulation} \ref{sec:example}

We consider $n =10$ points $x_i$ uniformly and randomly sampled in $[-10, 10] $. 
We synthesize also $Y$ data from a simple linear Regression :

$$Y =  \beta_1 X+\beta_0 + \epsilon $$

The values of the true parameters are $\beta_0= 7$ and $\beta_1= 4$  where we consider a random perturbation $\epsilon$ that we will take either Gaussian $N(0, \sigma^2 = 3)$  or uniform on $[-3 , 3]$ (so $\theta=3$, $\operatorname{var}(\epsilon)=\frac{\theta^2}{3} = 3 $). 
\noindent
\begin{figure}[H]
\centering
\includegraphics[width=150mm]{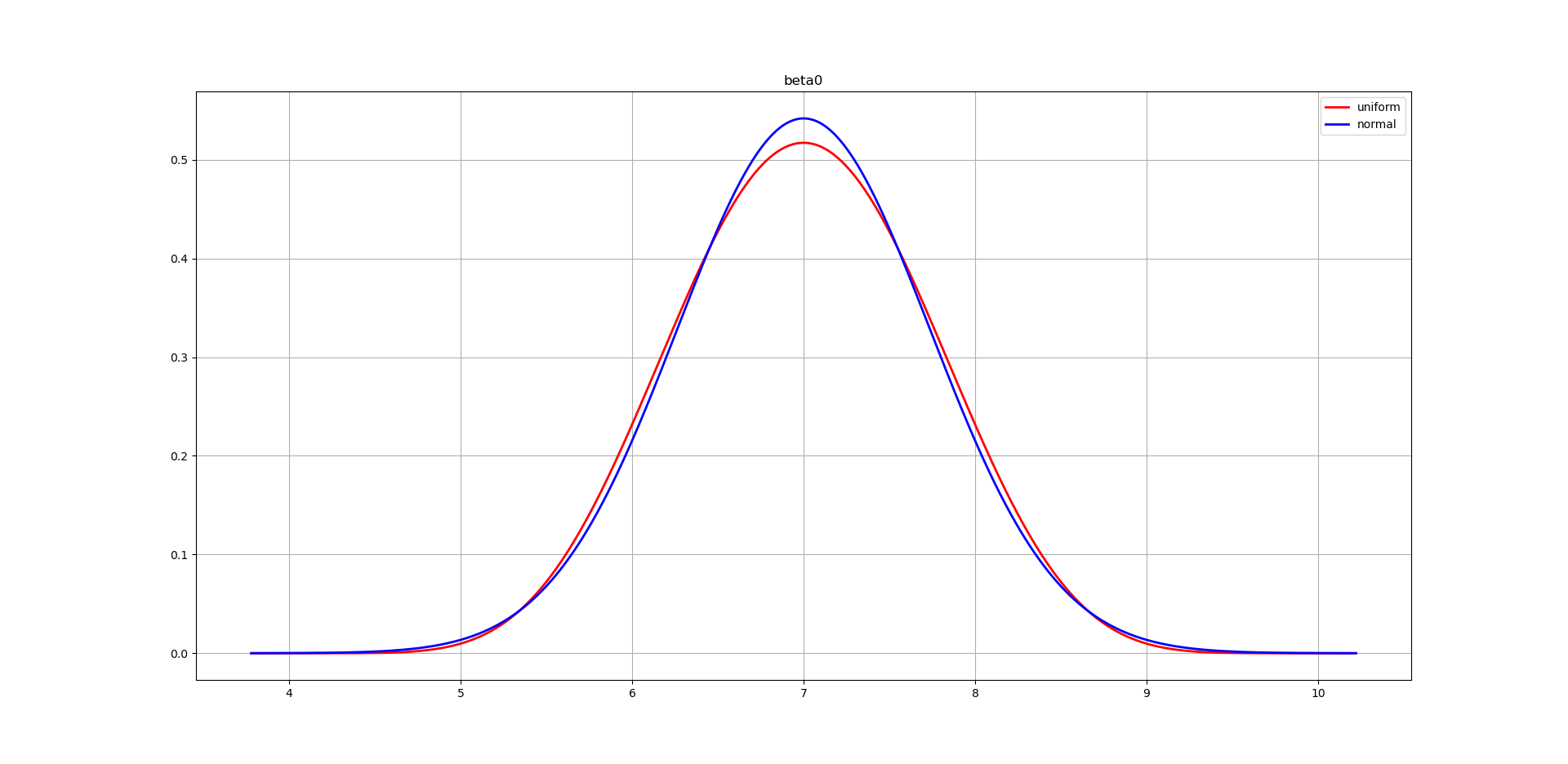}
\caption{Distributions for the estimator $\hat{\beta}_0$ }
\end{figure} 
\noindent
\begin{figure}[H]
\centering
\includegraphics[width=150mm]{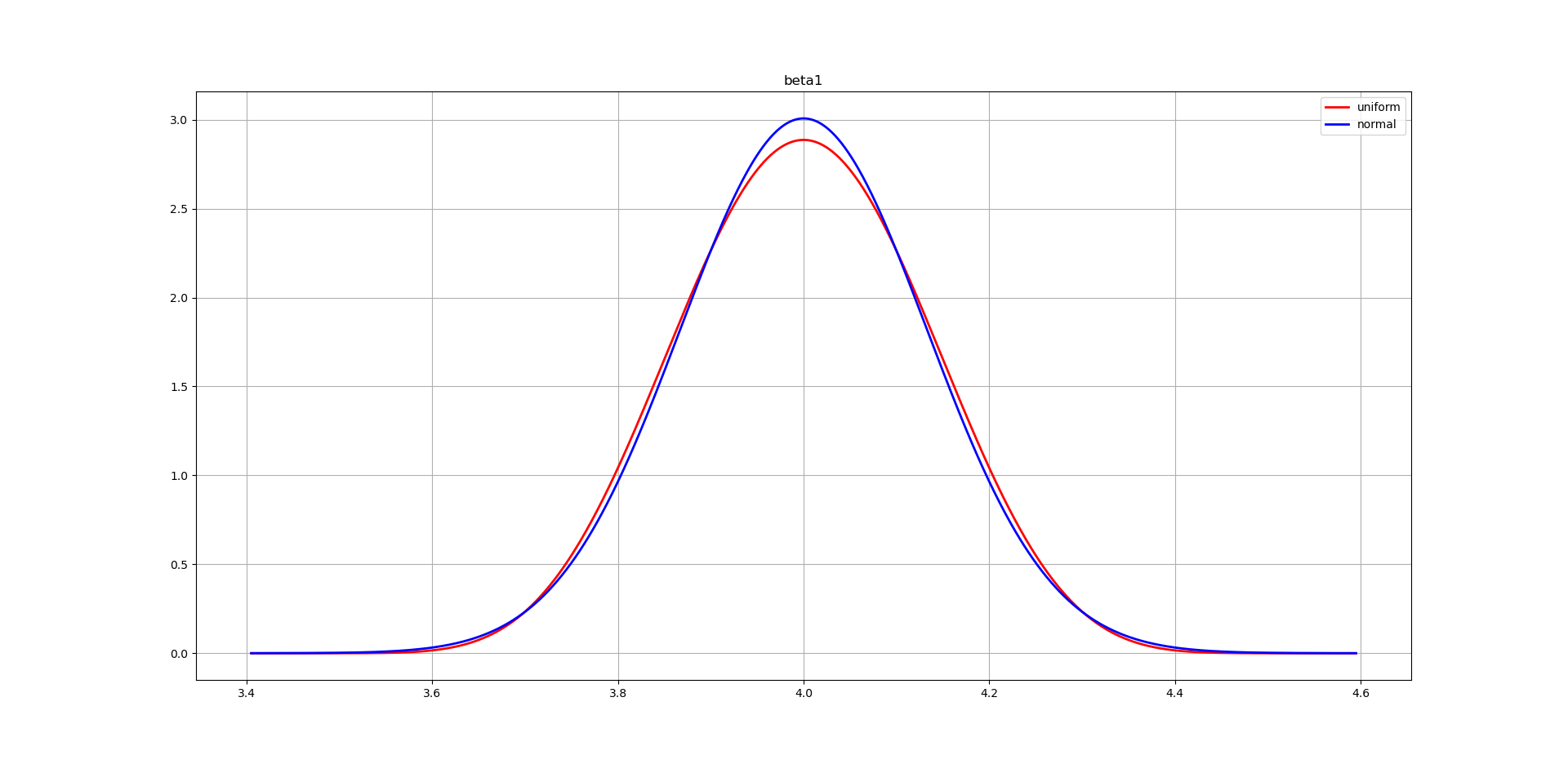}
\caption{Distributions for the estimator $\hat{\beta}_1$ }
\end{figure} 

The $95\%$ symmetric confidence intervals for $\hat{\beta_0}$ are of the form $[\beta_0-h, \beta_0+h]$. We show below the confidence intervals depending on $n$ in the case of the normal and uniform case, for  $\hat{\beta_0}$ and $\hat{\beta_1}$:

\begin{figure}[H]
\centering
\includegraphics[width=100mm]{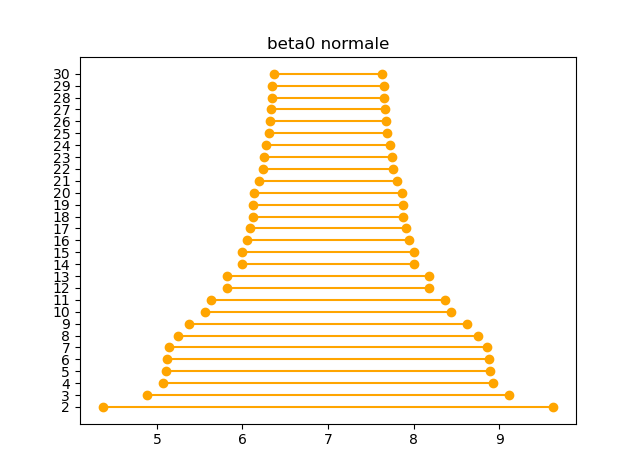}
\caption{$\beta_0$, $95\%$ confidence interval with normal error }
\end{figure} 
\noindent
\begin{figure}[H]
\centering
\includegraphics[width=100mm]{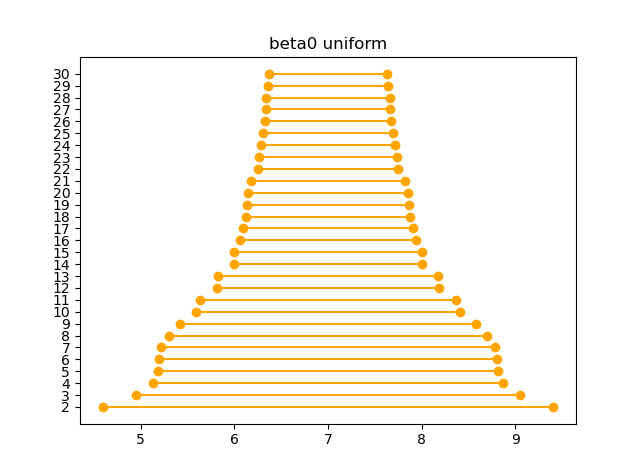}
\caption{$\beta_0$, $95\%$ confidence interval with uniform error}
\end{figure}
\newpage
\noindent
\begin{figure}[H]
\centering
\includegraphics[width=100mm]{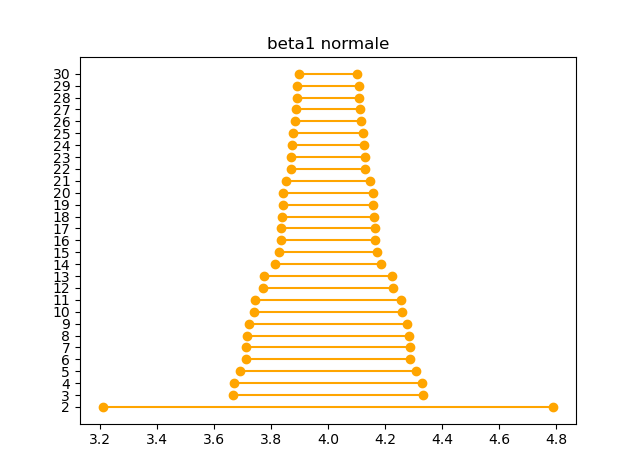}
\caption{ $\beta_1$, $95\%$ confidence interval with normal error }
\end{figure}
\noindent
\begin{figure}[H]
\centering
\includegraphics[width=100mm]{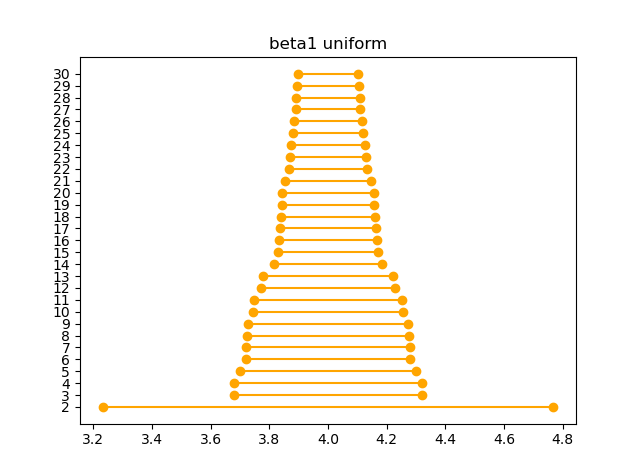}
\caption{$\beta_1$, $95\%$ confidence interval with uniform error}
\end{figure}

We see that the even for $n=10$, the distributions for the estimators are very close and that the confidence intervals are very close as well. However we note that we used the estimator from Theorem \ref{th-theta} to set $\sigma^2 = 3$. If we use the estimator assuming a Gaussian distribution, we get $\epsilon \approx N(0,1)$:
\noindent
\begin{figure}[H]
\centering
\includegraphics[width=100mm]{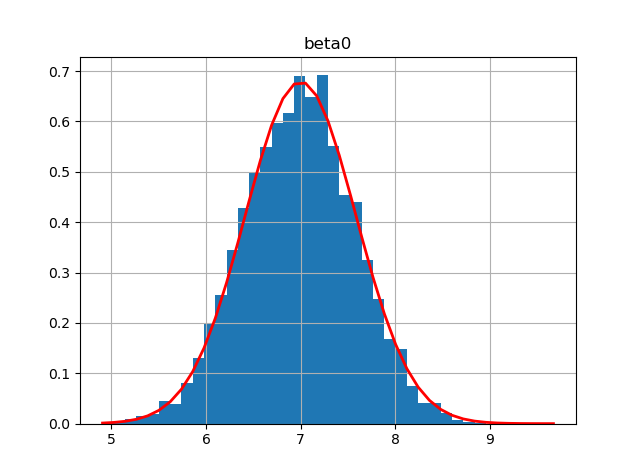}
\caption{$\beta_1$ distribution, assuming normal error $\sigma^2=1$.}
\end{figure} 
\begin{figure}[H]
\centering
\includegraphics[width=100mm]{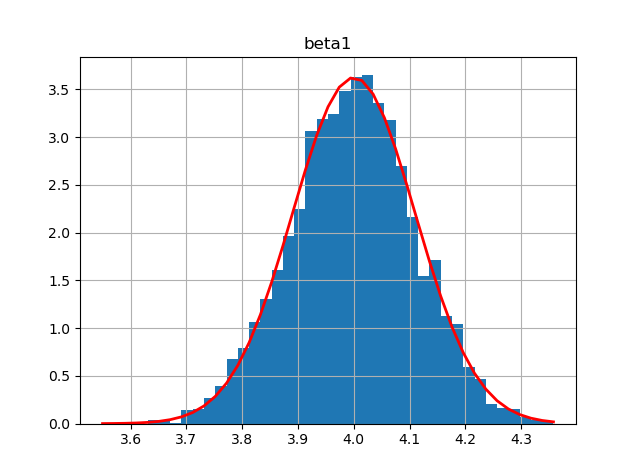}
\caption{$\beta_1$ distribution, assuming normal error $\sigma^2=1$.}
\end{figure}

We see here that the confidence intervals will differ a lot with respect to the uniform case, because the variance $\sigma^2$ is estimated assuming a Gaussian distribution with  $\frac{1}{n-2} ||Y-\hat{Y}||^2$.

\section{Conclusion}
We gave above the law for the least square estimator in the case of simple linear regression model with a uniform distribution on the perturbation term $\epsilon$. 
The uniform distribution hypothesis can be important to estimate the variance of the residual term. Visually it can be hard to assess if the model residuals are following a normal or a uniform distribution, however the variance parameter estimators differs significantly which leads to different confidence intervals for the parameters.
\bibliographystyle{plain}
\bibliography{References}

\hypertarget{bradl}{David M. Bradley and Ramesh C. Gupta. (2004). "On the distribution of the sum of n non-identically distributed uniform random variables" }.

\hypertarget{sadoo}{S. M. Sadooghi-Alvandi, A. R. Nematollahi and R. Habibi (2007). "On the distribution of the sum of independent uniform random variables"}. 

\end{document}